\documentclass{article}
\usepackage{amsfonts}
\usepackage{amsmath}

\input{tcilatex}
\begin{document}

\title{Free nearrings}
\author{Stefan Veldsman \\
Nelson Mandela University (South Africa)\\
and\\
La Trobe University (Australia)}
\maketitle

\begin{abstract}
An explicit construction of a free nearring as well as the free product (=
coproduct) of two nearrings are given in the variety of all, not necessarily
zero-symmetric, nearrings.
\end{abstract}

\begin{flushleft}
\textbf{AMS Subject Classification}: 16Y30

\textbf{Keywords: }nearring; free nearring
\end{flushleft}

\section{Introduction}

\noindent A \textit{nearring} $N$ is a quadruple $N=(N,+,\cdot ,0_{N})$
where $(N,+,0_{N})$ is a group with additive identity $0_{N},$ $(N,\cdot )$
is a semigroup and the multiplication distribute from the right over the
addition, i.e., $(a+b)c=ac+bc$ for all $a,b,c\in N.$ The group need not be
commutative. Strictly speaking, this is a \textit{right }nearring and there
is a corresponding notion of a \textit{left }nearring, but all our
considerations will be in the variety of right nearrings.

\bigskip

\noindent The prototype of all nearrings is the nearring $M(G),$ where $G$
is an additive group. $M(G)$ is the set of all functions $f:G\rightarrow G$
and the operations are pointwise addition and composition of functions. In a
nearring $N,$ $0a=0$ and $(-a)b=-(ab)$ for all $a,b\in N,$ but in general $%
a0 $ need not be $0$ and $a(-b)$ need not be $-(ab).$ A nearring $N$ is
called \textit{zero-symmetric} if $a0=0$ for all $a\in N.$ Any nearring $N$
can be decomposed as $N=N_{0}+N_{c}$ where $N_{0}:=\{a\in N\mid a0=0\}$ is
the \textit{zero-symmetric part} and $N_{c}:=\{a\in N\mid ab=a$ for all $%
b\in N\} $ is the \textit{constant part}. Note that $(N,\cdot ,0_{N})$ is
semigroup with left zero $0_{N}$ and $a0_{N}$ is also a left zero for all $%
a\in N.$ For more on nearrings, Pilz [7], Meldrum [6], Clay [2] or Ferrero
and Ferrero [3] can be consulted.

\bigskip

\noindent The construction of a free nearring as well as the free product of
two nearrings have been given by Meldrum [5] in the variety of
zero-symmetric left nearrings. Here we extend this work to the variety of
all, not necessarily zero-symmetric, nearrings. For this one could follow
the presentation in [5] almost verbatim. But note that here we work in the
variety of all \textit{right nearrings} and provide more details to, amongst
others, take care of the presence of constant elements.

\bigskip

\section{The free nearring}

\noindent Let $X$ be a non-empty set. By $F(X)$ we denote the free nearring
on the set $X.$

\bigskip

\noindent \textbf{2.1}\qquad Let $(R,\cdot ,0_{R})$ be the free
(multiplicative) semigroup with left zero $0_{R}$ over $X$ in the category
of semigroups with left zero. The objects in this category are $(S,\cdot
,0_{S})$ where $S$ is a semigroup with respect to the binary operation $%
\cdot $ and $0_{S}\in S$ with $0_{S}a=0_{S}$ for all $a\in S.$ A morphism is
given by $f:(S_{1},\cdot _{1},0_{1})\rightarrow (S_{2},\cdot _{2},0_{2})$
which is a semigroup homomorphism with $f(0_{1})=0_{2}.$

\noindent This free semigroup with left zero $(R,\cdot ,0_{R})$ has an
explicit description as: $R=R_{w}\cup R_{0}$ where

\qquad $\qquad R_{w}:=\{$words $x_{1}x_{2}...x_{n}\mid n\geq 1,x_{i}\in X\}$

\noindent is the free semigroup over $X$ with concatenation as operation and

\qquad \qquad $R_{0}:=\{0_{R}\}\cup \{w0_{R}\mid w\in R_{w}\}$

\noindent for some symbol $0_{R}$ (not in $X)$ and formal products $%
w0_{R},w\in R_{w}.$ For $u,v\in R,$ the product $uv$ in $R$ is given by:

\qquad $uv:=\left\{ 
\begin{tabular}{l}
$uv$ product in $R_{w}$ for $u,v\in R_{w}$ \\ 
$u0_{R}$ formal product in $R_{0}$ for $u\in R_{w},v=0_{R}$ \\ 
$(uw)0_{R}$ formal product in $R_{0}$ for $u\in R_{w},v=w0_{R}$ \\ 
$u$ for $u\in R_{0},v\in R$%
\end{tabular}%
\right. .$

\bigskip

\noindent \textbf{2.2}\qquad Define a partial addition $+$ on $R$ by $%
u+0_{R}=u=0_{R}+u$ for all $u\in R.$

\noindent On the structure $(R,+,\cdot ,0_{R}),$ we will now apply the
inductive steps as given by Meldrum to get $F(X).$ These steps are:

\noindent Let

\begin{tabular}{l}
$A_{0}=R;$ \\ 
$B_{0}=Fg(R,+)$ the free additive group on $(R,+)$ extending the addition on 
$R;$ \\ 
$C_{0}=B_{0}-A_{0}$; and \\ 
$D_{0}=R\cdot C_{0}$ where $R\cdot C_{0}$ indicates the formal product $%
r\cdot c$ of $r\in R$ and $c\in C_{0}.$%
\end{tabular}

\bigskip

\noindent If $A_{n-1},B_{n-1},C_{n-1}$ and $D_{n-1}$ have been defined for $%
n\geq 1,$ let

\bigskip

\begin{tabular}{l}
$A_{n}=B_{n-1}\cup D_{n-1};$ \\ 
$B_{n}=B_{n-1}\ast Fg(D_{n-1})$ where $Fg(D_{n-1})$ is the free additive
group on $D_{n-1}$ \\ 
\ \ and $\ast $ denotes the free product of the two groups; \\ 
$C_{n}=B_{n}-A_{n};$ and \\ 
$D_{n}=R\cdot C_{n}$ where $R\cdot C_{n}$ indicates the formal product $%
r\cdot c$ of $r\in R$ and $c\in C_{n}.$%
\end{tabular}

\bigskip

\noindent Let $F(X):=\tbigcup\limits_{n=0}^{\infty }B_{n}$. Every $B_{n-1}$
is naturally embedded in $B_{n}.$ Hence $F(X)$ is the union of a tower of
groups which again is a group; it is actually the free product of the group $%
B_{0}$ and a free group. Meldrum then defines the multiplication in $F(X)$
and shows that it is right distributive over the addition and also
associative. Below we will follow these steps and check the required
properties.

\noindent We write the additive identity in $F(X)$ as $0.$ Note that $%
B_{n}=C_{n}\cup A_{n}$ and $B_{n}-B_{n-1}=C_{n}\cup D_{n-1}.$ If we let $%
D_{-1}:=R,$ then $\tbigcup\limits_{n=-1}^{\infty }D_{n}$ is a generating set
for $F(X)=\tbigcup\limits_{n=0}^{\infty }B_{n}$ since

\qquad $B_{0}=Fg(D_{-1},+),$

$\qquad B_{1}=B_{0}\ast Fg(D_{0})=Fg(D_{-1},+)\ast Fg(D_{0}),$

$\qquad B_{2}=B_{1}\ast Fg(D_{1})=Fg(D_{-1},+)\ast Fg(D_{0})\ast Fg(D_{1}),$
etc.

\bigskip

\noindent \textbf{2.3}\qquad Firstly we describe the two free constructions
required in the inductive steps:

\noindent (a) $B_{0}=Fg(R,+)$ the free additive group on $(R,+)$ extending
the addition on $R$: Let $R^{-}:=\{-r\mid r\in R-\{0_{R}\}\}$ where each $-r$
is just a symbol not in $R.$ Let $w=r_{1}+r_{2}+...+r_{n}$ be a word over $%
R\cup R^{-1},$ written additively, where $n\geq 1$ and $r_{i}\in R\cup
R^{-1} $. The word $w$ is reduced if:

(i) If $r_{i}\in \{b,-b\}$ for some $b\in R-\{0_{R}\}$ with $r_{i+1}\in
\{b,-b\}-\{r_{i}\},$ replace $r_{i}+r_{i+1}$ with $0_{R};$

(ii) If $r_{i}=0_{R},$ replace $r_{i-1}+0_{R}$ with $r_{i-1}$ and $%
0_{R}+r_{i+1}$ with $r_{i+1}$ (of course, ignore the replacement if $r_{i-1}$
or $r_{i+1}$ does not make sense).

\noindent Then $B_{0}$ is the set of all reduced words over $R\cup R^{-}$
and with the operation of additive concatenation of two reduced words
followed by reduction if necessary, it is a group. The additive identity in $%
B_{0}$ is $0_{R}$ which we write as $0.$

\bigskip

\noindent (b) $B_{n}\ast Fg(D_{n})$ where $Fg(D_{n})$ is the free additive
group on $D_{n}$ and $\ast $ denotes the free product of the two groups.

\noindent Firstly, $Fg(D_{n})$ consists of all reduced words $%
w=d_{1}+d_{2}+...+d_{n},n\geq 0,$ over $D_{n}\cup D_{n}^{-}$ where $%
D_{n}^{-} $ is the set of symbols $D_{n}^{-}=\{-d\mid d\in D_{n}\}$ and the
reduction is:

If $d_{i}\in \{d,-d\}$ for some $d\in D_{n}$ with $d_{i+1}\in
\{d,-d\}-\{d_{i}\},$ replace $d_{i}+d_{i+1}$ with $0_{D_{n}}.$

\noindent Here the empty word is denoted by $0_{D_{n}}$ and it is the
additive identity for the operation of concatenation followed by reduction
if necessary.

\noindent Secondly, for the coproduct $B_{n}\ast Fg(D_{n}):$ Consider all
words $u=u_{1}+u_{2}+...+u_{n},n\geq 0,$ over $B_{n}\cup D_{n}\cup D_{n}^{-}$
subject to the reduction:

(i) If $u_{i}\in \{d,-d\}$ for some $d\in D_{n}$ with $u_{i+1}\in
\{d,-d\}-\{u_{i}\},$ replace $u_{i}+u_{i+1}$ with $0_{D_{n}};$

(ii) If $u_{i}+u_{i+1}$ is already defined (eg., both $u_{i},u_{i+1}$ are in
the group $B_{n}$), say $u_{i}+u_{i+1}=u_{i}^{^{\prime }},$ replace $%
u_{i}+u_{i+1}$ in $u$ with $u_{i}^{^{\prime }}.$

(iii) Remove all occurances of $0_{B_{n}}$ and $0_{D_{n}}$ from $u.$

\noindent Then $B_{n}\ast Fg(D_{n})$ consists of all the reduced words,
including the empty word denoted by $0$, with group operation the
concatenation of words followed by reduction and $0$ is the additive
identity.

\bigskip

\noindent \textbf{2.4}\qquad Next it is shown how the product in $F(X)$ is
defined. Note that for any $b\in F(X),$ we often consider the three
possibilites $b\in R\cup C_{n}\cup D_{n}$ for some $n\geq 0.$ Indeed, for $%
b\in F(X)=\tbigcup\limits_{n=0}^{\infty }B_{n},$ if $b\in B_{0},$ then $b\in
R$ or $b\in B_{0}-R=C_{0}.$ Otherwise $b\in B_{n}-B_{n-1}=C_{n}\cup D_{n-1}$
for some $n\geq 1.$

\noindent Step 1. Define the product $RF(X).$ This is done in several
stages. Let $a\in R,b\in F(X).$

\noindent (a) If $a=0,$ then $ab=0b:=0$ for all $b\in F(X).$

\noindent Suppose thus $a\neq 0$ in what follows.

\noindent (b) If $b\in R,$ then $ab$ is just the product in $R.$

\noindent (c) If $b\in C_{n}$ for some $n\geq 0,$ then $ab:=a\cdot b$ the
formal product in $R\cdot C_{n}=D_{n}.$

\noindent (d) If $b\in D_{n}=R\cdot C_{n}$ for some $n\geq 0,$ then $%
b=r\cdot c$ for some $r\in R$ and $c\in C_{n}.$ Then $ab=a(r\cdot
c):=(ar)\cdot c$ the formal product of $ar\in R$ and $c\in C_{n}.$

\noindent Thus $ab$ is defined for all $a\in R,b\in F(X).$

\bigskip

\noindent Step 2. The product $R^{-}F(X)$: Let $a\in R^{-},$ say $a=-r,r\in
R $, and $b\in F(X).$ Then $ab=(-r)b:=-(rb)$ which is well-defined with $%
rb\in RF(X)\subseteq F(X)$ by Step 1 and the additive inverse of $rb$ in the
group $F(X).$ This means, for $d\in D_{n}$ and its additive inverse $-d\in
D_{n}^{-}$ in the group $Fg(D_{n}),$ and hence in $F(X),$ we have $d=r\cdot
c $ for some $r\in R,c\in C_{n}$ and so $-d=-(r\cdot c)=-(rc)=(-r)c$ by Step
1(c) and the definition here in Step 2. Hence $D_{n}^{-}=\{-(rc)\mid r\in
R,c\in C_{n}\}=\{(-r)c\mid r\in R,c\in C_{n}\}.$

\bigskip

\noindent Step 3. The product $B_{0}F(X)$: Let $a\in B_{0}=Fg(R,+),b\in
F(X). $ Then $a=a_{1}+a_{2}+...+a_{n}$ where $a_{i}\in R\cup R^{-}.$ By
Steps 1 and 2, $a_{i}b$ is defined in $F(X)$ for all $i,$ and so is the sum $%
a_{1}b+a_{2}b+...+a_{n}b$ in the group $F(X).$ Let $%
ab=(a_{1}+a_{2}+...+a_{n})b:=a_{1}b+a_{2}b+...+a_{n}b.$ Thus $B_{0}F(X)$ is
defined.

\bigskip

\noindent Step 4. Inductive step to define $B_{n}F(X)$ assuming that $%
B_{n-1}F(X),$ has been defined for some $n\geq 1$: Let $b\in F(X).$

\noindent Let $a\in B_{n}-B_{n-1}=C_{n}\cup D_{n-1}.$ There are two cases to
consider.

\noindent (a) If $a\in D_{n-1}=R\cdot C_{n-1},$ then $a=r\cdot c$ for some $%
r\in R,c\in C_{n-1}\subseteq B_{n-1}.$ Let $ab=(r\cdot c)b:=r(cb)$ which is
well-defined by Step 1 since $cb\in B_{n-1}F(X)\subseteq F(X)$ by the
induction assumption.

\noindent (b) If $a\in C_{n}\subseteq B_{n}=B_{n-1}\ast Fg(D_{n-1}),$ then $%
a=a_{1}+a_{2}+...+a_{n}$ where $a_{i}\in B_{n-1}\cup D_{n-1}\cup
D_{n-1}^{-}. $ We consider the three cases.

(i) If $a_{i}\in B_{n-1},$ then $a_{i}b$ is defined by the induction
assumption.

(ii) If $a_{i}\in D_{n-1},$ then $a_{i}b$ is already defined as in (a) above.

(iii) If $a_{i}\in D_{n-1}^{-}$, say $a_{i}=-d$ for some $d\in D_{n-1},$
then $a_{i}b=(-d)b:=-(db)$ is well-defined by (a) and as an additive inverse
in the group $F(X).$

\bigskip

\noindent \textbf{2.5}\qquad By 2.4, the group $F(X)$ has a well-defined
product. This product is distributive from the right over the addition as is
shown next. If $a,b,c\in F(X)=\tbigcup\limits_{n=0}^{\infty }B_{n},$ we may
choose $n\geq 0$ with $a,b\in B_{n}$. We know $B_{n}$ is a free group;
either $Fg(R,+)$ or $B_{n-1}\ast Fg(D_{n-1}).$ In both cases, the products $%
ac$ and $bc$ were defined by distributivity from the right over the sum of
the generators for $a$ and $b$, hence we should have $(a+b)c=ac+bc.$
However, there could be some concern if the sum $a+b$ is subject to
reductions. But this is not a problem. For example, suppose $%
a=a_{1}+a_{2}+...+a_{n}$ and $b=b_{1}+b_{2}+...+b_{m}$ where $a_{i},b_{j}\in
R\cup R^{-}$ for all $i$ and $j,$ or $a_{i},b_{j}\in B_{n-1}\cup D_{n-1}\cup
D_{n-1}^{-}$ for all $i$ and $j,$ with a reduction $a_{n}+b_{1}=d$ for some $%
d,$ which could be $0.$ Thus

\qquad $(a+b)c=(a_{1}+a_{2}+...+a_{n-1}+d+b_{2}+...+b_{m})c$

$\qquad \qquad =a_{1}c+a_{2}c+...+a_{n-1}c+dc+b_{2}c+...+b_{m}c$ while

$\qquad ac+bc=a_{1}c+a_{2}c+...+a_{n-1}c+a_{n}c+b_{1}c+b_{2}c+...+b_{m}c$

$\qquad \qquad =(a+b)c$ since $a_{n}c+b_{1}c=dc$ by definition of the
product $dc.$

\noindent If $a_{n-1}+d,$ $d+b_{2}$ or when $d=0,a_{n-1}+b_{2}$ can be
reduced, we may repeat this argument as many times as is necessary to
conclude that the product distributes over the addition in $B_{n}.$

\noindent Thus the multiplication in $F(X)$ is right distributive over the
addition in $F(X).$

\bigskip

\noindent \textbf{2.6}\qquad By the distribitivity, for any $a,b\in F(X),$ $%
0=0b=(a+(-a))b=ab+(-a)b$ and hence $(-a)b=-(ab).$

\bigskip

\noindent \textbf{2.7}\qquad The final step to validate that $F(X)$ is a
nearring, is to show that the product is associative. By a result of Laxton
and Lockheart [3], it is sufficient to show that $a(bc)=(ab)c$ for all $%
b,c\in F(X)$ and for all $a$ from a generating set for $F(X).$ For this we
use $\tbigcup\limits_{n=-1}^{\infty }D_{n}$ and we proceed by induction on $%
n $.

\bigskip

\noindent Step 1. Let $a\in D_{-1}=R.$ We show $a(bc)=(ab)c$ for all $b,c\in
F(X)$ and in doing this, we consider the three possibilites for $b\in F(X),$
namely $b\in R,b\in C_{k}$ or $b\in D_{k}$ for some $k\geq 0$ (refer the
opening lines of 2.4 above).

\noindent (a) For $b\in R,$ we consider the three subcases:

(i) $c\in R:$ $a(bc)=(ab)c$ in the semigroup $(R,\cdot );$

(ii) $c\in C_{t}$ for some $t\geq 0:$ $a(bc)=a(b\cdot c)=(ab)\cdot c=(ab)c$
by definition in 2.4, Step 1(d).

(iii) $c\in D_{t}=R\cdot C_{t}$ for some $t\geq 0,$ say $c=r\cdot c^{\prime
} $ for some $r\in R,c^{\prime }\in C_{t}:$ $a(bc)=a(b(r\cdot c^{\prime
}))=a((br)\cdot c^{\prime })=(a(br))\cdot c^{\prime }=((ab)r)\cdot c^{\prime
}=(ab)(r\cdot c^{\prime })=(ab)c$ where we have used 2.4 Step 1(d) several
times as well as the associativity in the semigroup $(R,\cdot ).$

\noindent (b) For $b\in C_{k}:(ab)c=(a\cdot b)c=a(bc)$ by 2.4 Step1(c) and
Step 4.

\noindent (c) For $b\in D_{k}=R\cdot C_{k},$ say $b=r\cdot c^{\prime }$ for
some $r\in R,c^{\prime }\in C_{k}\subseteq B_{k}$ we have:

\qquad \qquad 
\begin{tabular}{lll}
$(ab)c$ & $=$ & $(a(r\cdot c^{\prime }))c$ \\ 
& $=$ & $((ar\cdot c^{\prime })c$ by 2.4 Step 1(d) \\ 
& $=$ & $(ar)(c^{\prime }c)$ by 2.4 Step 4 \\ 
& $=$ & $a(r(c^{\prime }c))$ by 2.4 Step 4 \\ 
& $=$ & $a((r\cdot c^{\prime })c))$ by 2.4 Step 4 \\ 
& $=$ & $a(bc).$%
\end{tabular}

\bigskip

\noindent Step 2. Let $n\geq 1$ and suppose for all $m$ with $0\leq m<n$ the
equality $x(yz)=(xy)z$ holds for all $x\in D_{m},y,z\in F(X).$

\noindent Let $a\in D_{n},b,c\in F(X).$ Then $a=v\cdot u$ for some $v\in R$
and $u\in C_{n}\subseteq B_{n},$ hence $u=\sum u_{i}$ for $u_{i}\in
B_{n-1}\cup D_{n-1}\cup D_{n-1}^{-}.$ Therefor,

\qquad \qquad 
\begin{tabular}{lll}
$(ab)c$ & $=$ & $((v\cdot u)b)c$ \\ 
& $=$ & $((v(ub))c$ by 2.4 Step 5 \\ 
& $=$ & $(v((\sum u_{i})b))c$ \\ 
& $=$ & $(v(\sum (u_{i}b))c$ by 2.4 Step 4 \\ 
& $=$ & $v(\sum (u_{i}b)c)$ by 2.4 Step 4%
\end{tabular}

\noindent and

\qquad \qquad 
\begin{tabular}{lll}
$a(bc)$ & $=$ & $(vu)(bc)$ by 2.4 Step 1(c) \\ 
& $=$ & $v(u(bc))$ by 2.4 Step 4 \\ 
& $=$ & $v((\sum u_{i})(bc))$ \\ 
& $=$ & $v(\sum u_{i}(bc))$ by 2.4 Step 4.%
\end{tabular}

\noindent The equality $(ab)c=a(bc)$ will follow from the equality $%
(u_{i}b)c=u_{i}(bc)$ for $u_{i}\in B_{n-1}\cup D_{n-1}\cup D_{n-1}^{-}$
which we will now do. This is done in turn for the three possible choices of 
$u_{i}:$

\noindent (a) For $u_{i}\in D_{n-1},$ $(u_{i}b)c=u_{i}(bc)$ follows from the
induction assumption.

\noindent (b) For $u_{i}\in D_{n-1}^{-},$ say $u_{i}=-d,d\in D_{n-1},$ we
have by 2.6 and (a),

$\qquad (u_{i}b)c=((-d)b)c=(-(db))c=-((db)c)=-(d(bc))=(-d)(bc)=u_{i}(bc).$

\noindent (c) For $u_{i}\in B_{n-1}:$ We know that $B_{n-1}$ is generated as
a group by $\tbigcup\limits_{m=-1}^{n-1}D_{m}$ which means that $u_{i}$ can
be expressed in terms of elements $d_{j}\in D_{m_{j}}$ for $0\leq m_{j}\leq
n-1.$ By the induction assumption, we know $d_{j}(bc)=(d_{j}b)c$ for all
these $d_{j}$'s. Thus $(u_{i}b)c=u_{i}(bc)$ follows.

\bigskip

\noindent \textbf{2.8}\qquad We may thus conclude that the product is
associative and that $F(X)$ is a nearring. Lastly we show $F(X)\ $is the
free nearring on $X.$ Let $\iota :X\rightarrow F(X)$ be the inclusion
function and let $g:X\rightarrow N$ be any function from $X\ $to a nearring $%
N.$ We show that $g$ has a unique extension to a nearring homomorphism $%
f:F(X)\rightarrow N$ (so $f\circ \iota =g):$

\noindent Step 1. Define $f(a)$ for $a\in B_{0}:$

\noindent If $a\in R,$ then $a=x_{1}x_{2}...x_{k}$ or $a=0_{R}$ or $%
a=x_{1}x_{2}...x_{k}0_{R}$ for $k\geq 1$ and $x_{i}\in X.$ Let

$f(x_{1}x_{2}...x_{k}):=g(x_{1})g(x_{2})...g(x_{k});$

$f(0_{R})=f(0):=0_{N};$ and

$f(x_{1}x_{2}...x_{k}0_{R}):=g(x_{1})g(x_{2})...g(x_{k})0_{N}.$

\noindent With $f(a)\in N$ defined for all $a\in R,$ we let

$f(-a):=-f(a),$ the additive inverse of $f(a)$ in the nearring $N.$

\noindent For $a\in B_{0}-R,$ say $a=a_{1}+a_{2}+...+a_{k}$ where $k\geq 1$
and $a_{i}\in R\cup R^{-}$, let

$\qquad f(a):=f(a_{1})+f(a_{2})+...+f(a_{k}).$

\bigskip

\noindent Step 2. Define $f(a)$ for $a\in B_{n}$ and $n\geq 1$ subject to
the assumption that $f(b)$ is defined for all $b\in B_{n-1}:$ Let $a\in
B_{n}-B_{n-1}=C_{n}\cup D_{n-1}.$

If $a\in D_{n-1}=R\cdot C_{n-1},$ say $a=r\cdot c$ where $r\in R$ and $c\in
C_{n-1}\subseteq B_{n-1},$ let $f(a)=f(rc):=f(r)f(c)$ which is well-defined
as the product in the nearring $N$ of $f(r)\in N$ and $f(c)\in N$ by Step 1
and the induction assumption. Let $f(-a):=-f(a).$

If $a\in C_{n}\subseteq B_{n},$ then $a=a_{1}+a_{2}+...+a_{k}$ where $k\geq
1 $ and $a_{i}\in B_{n-1}\cup D_{n-1}\cup D_{n-1}^{-}.$ Since $f(a_{i})$ is
defined for all $i$ by the induction assumption and the first part of this
step, we let $f(a):=f(a_{1})+f(a_{2})+...+f(a_{k}).$

\bigskip

\noindent Step 3. By Steps 1 and 2, it follows that $f:F(X)\rightarrow N$ is
a well-defined function with $f\circ \iota =g.$ As can be seen from its
definition above, $f$ is a nearring homomorphism, taking care to check
possible reductions in the sum $a+b$ and product $ab$ when showing $%
f(a+b)=f(a)+f(b)$ and $f(ab)=f(a)f(b)$ respectively. It is also
straightforward to show that any nearring homomorphism $h:f(X)\rightarrow N$
which coincides with $g$ on $X,$ must coincides with $f$ on $F(X).$

\noindent Thus $F(X)$ is the free nearring on the set $X.$

\bigskip

\noindent This nearring does not have an identity. For the free nearring
with an identity on the set $X,$ we can proceed as in the next example as
suggested by Meldrum [5].

\bigskip

\noindent \textbf{Example.} Let $X$ be a non-empty set. Let $(R,\cdot
,1_{R},0_{R})$ be the free multiplicative semigroup with identity $1_{R}$
and left zero $0_{R}$ over $X$ in the category of semigroups with identity
and left zero. The objects in this category are $(S,\cdot ,1_{S},0_{S})$
where $S$ is a semigroup with respect to the binary operation $\cdot $, $%
1_{S}\in S$ with $a1_{S}=a=a1_{S}$ for all $a\in S$ and $0_{S}\in S$ with $%
0_{S}a=0_{S}$ for all $a\in S.$ A morphism is given by $f:(S_{1},\cdot
_{1},1_{1},0_{1})\rightarrow (S_{2},\cdot _{2},1_{2},0_{2})$ which is a
semigroup homomorphism with $f(1_{1})=1_{2}$ and $f(0_{1})=0_{2}.$

\noindent This free semigroup with identity and left zero $(R,\cdot
,1_{R},0_{R})$ has an explicit description as: $R=R_{w}\cup R_{0}\cup
\{1_{R}\}$where

\qquad $\qquad R_{w}:=\{$words $x_{1}x_{2}...x_{n}\mid n\geq 1,x_{i}\in X\}$

\noindent is the free semigroup over $X$ with concatenation as operation;

\qquad \qquad $R_{0}:=\{0_{R}\}\cup \{w0_{R}\mid w\in R_{w}\}$

\noindent for some symbol $0_{R}$ (not in $X)$ and formal products $%
w0_{R},w\in R_{w};$ and

$\qquad \qquad 1_{R}$ is a symbol different to $0_{R}$ and not in $X.$

\noindent For $u,v\in R,$ the product $uv$ in $R$ is given by:

\qquad $uv:=\left\{ 
\begin{tabular}{l}
$u$ if $v=1_{R}$ \\ 
$v$ if $u=1_{R}$ \\ 
$uv$ product in $R_{w}$ for $u,v\in R_{w}$ \\ 
$u0_{R}$ formal product in $R_{0}$ for $u\in R_{w},v=0_{R}$ \\ 
$(uw)0_{R}$ formal product in $R_{0}$ for $u\in R_{w},v=w0_{R}$ \\ 
$u$ for $u\in R_{0},v\in R$%
\end{tabular}%
\right. .$\bigskip

\noindent With this $R,$ write $1$ and $0$ for $1_{R}$ and $0_{R}$
respectively and then proceed with 2.2 and the subsequent steps to construct 
$F(X),$ subject to changing the definition of $D_{n}$ from $D_{n}=R\cdot
C_{n}$ to $D_{n}=(R-\{1\})\cdot C_{n}$ and using $1$ as the multiplicative
identity throughout the process.

\bigskip

\section{The free product of two nearrings}

\noindent Here the free product of two nearrings $A$ and $B$ will be
described. Again we follow Meldrum [4] who uses the same steps as in the
construction above of the free nearring. The starting point $R$ is in this
case different, being the free product of two semigroups. Remember, Meldrum
works with zero-symmetric left nearrings but here we work in the variety of
all, not necessarily zero-symmetric, right nearrings.

\bigskip

\noindent Let $A=(A,+,\cdot ,0_{A})$ and $B=(B,+,\cdot ,0_{B})$ be two
nearrings, right distributive and not necessarily zero-symmetric. Without
loss of generality, we may assume $A$ and $B$ are disjoint. Let $R=(R,\cdot
,0_{R})$ be the free product of the two semigroups with left zero $(A,\cdot
,0_{A})$ and $(B,\cdot ,0_{B}).$ An explicit construction of $R$ can be
given as follows. A word $w=x_{1}x_{2}...x_{n}$ over $A\cup B,$ $x_{i}\in
A\cup B,$ is reduced if:

(i) If $x_{i},x_{i+1}\in A$ $($respt. $B)$ and $x_{i}x_{i+1}=x_{i}^{\prime
}\in A$ $($respt. $B),$ replace $x_{i}x_{i+1}$ with $x_{i}^{\prime }$ in $w.$

(ii) If $x_{i}\in \{0_{A},0_{B}\},$ replace $x_{1}x_{2}...x_{n}$ with $%
x_{1}x_{2}...x_{i-1}0_{A}.$

\noindent Let $R$ be the set of all reduced words over $A\cup B$ with
operation concatenation followed by reduction. Then $R$ is a semigroup with
left zero $0_{R}:=0_{A}(=0_{B}).$ It can be checked that $R$ is the
coproduct of $(A,\cdot ,0_{A})$ and $(B,\cdot ,0_{B})$ in the category of
semigroups with left zero as described in 2.1 above.

\bigskip

\noindent Since $A\cup B\subseteq R,$ there is a partial addition $+$
defined on $R$ given by:

For $a_{1},a_{2}\in A,$ $a_{1}+a_{2}\in A\subseteq R$ is just the addition
in the nearring $A;$

for $b_{1},b_{2}\in B,$ $b_{1}+b_{2}\in B\subseteq R$ is just the addition
in the nearring $B;$ and,

by definition, for all $r\in R,r+0_{R}=r=0_{R}+r.$

\noindent Thus $R=(R,+,\cdot ,0_{R})$ is an algebraic structure with $%
(R,\cdot ,0_{R})$ a semigroup with left zero $0_{R},(R,+,0_{R})$ has a
partial addition $+$ on $R$ with additive identity $0_{R},$ and $R$ contains
the two nearrings $A=(A,+,\cdot ,0_{R})$ and $B=(B,+,\cdot ,0_{R}).$ For $%
w\in R,$ we define

\qquad $-w=\left\{ 
\begin{array}{l}
-a\text{ additive inverse in }A\text{ if }w=a\in A,a\neq 0_{R} \\ 
-b\text{ additive inverse in }B\text{ if }w=b\in B,b\neq 0_{R} \\ 
0_{R}\text{ if }w=0_{R} \\ 
-w\text{ formal symbol if }w\in R-A\cup B%
\end{array}%
\right. .$

\noindent Let $R^{-}:=\{-w\mid w\in R-\{0_{R}\}\}.$

\noindent We can then construct the free nearring $F(R)$ on $R$ following
the Meldrum inductive steps as given in 2.2 with $A_{0}=R$ and $%
B_{0}=Fg(R,+),$ the latter the free additive group on $R$ extending the
partial addition $+$ to all of $R.$ Although this free group was already
described in 2.3, the one here is for a different $R$ and it may be
instructive to explicitly repeat it here for this case. A word $w$ over $%
R\cup R^{-1}$ is a sum $w=r_{1}+r_{2}+...+r_{n}$ with $r_{i}\in R\cup R^{-1}$%
, $n\geq 1.$ The word $w$ is reduced if:

(i) If $r_{i}+r_{i+1}$ is already defined in $R,$ say $r_{i}+r_{i+1}=r_{i}^{%
\prime },$ replace $r_{i}+r_{i+1}$ with $r_{i}^{\prime }$ in $w.$

(ii) If $r_{i}\in \{u,-u\}$ for some $u\in R$ with $r_{i+1}\in
\{u,-u\}-\{r_{i}\},$ replace $r_{i}+r_{i+1}$ with $0_{R}.$

\noindent Then $B_{0}$ is the set of all reduced words over $R\cup R^{-1}$.
It is a group with respect to the operation of (additive) concatenation of
two reduced words followed by reduction if necessary. The additive identity
in $B_{0}$ is $0_{R}$ which is written as $0$ in the sequel.

\bigskip

\noindent As in the previous section, it follows that $F(R)$ is a free
additive group endowed with a well-defined associative multiplication which
is right distributive over the addition. We show that $F(R)\ $is the
coproduct of the two nearrings $(A,+,\cdot ,0_{A})$ and $(B,+,\cdot ,0_{B})$
in the category of nearrings.

\noindent The maps $\iota _{A}:A\rightarrow F(R)$ and $\iota
_{B}:B\rightarrow F(R)$ defined by $\iota _{A}(a)=a$ for all $a\in A$ and $%
\iota _{B}(b)=b$ for all $b\in B$ are injective nearring homomorphisms. Let $%
N$ be a nearring and let $f_{A}:A\rightarrow N$ and $f_{B}:B\rightarrow N$
be two homomorphisms. We know $\iota _{A}(0_{A})=0_{A}=0=0_{B}=\iota
_{B}(0_{B})$ and $f_{A}(0_{A})=0_{N}=f_{B}(0_{B}).$ Define a map $%
f:F(R)\rightarrow N$ inductively by the following steps:

\noindent (a) Let $x\in F(R).$ If $x=0,$ let $f(x):=0_{N}.$ Suppose thus $%
x\neq 0.$

If $x\in A,$ say $x=a\in A,$ let $f(a):=f_{A}(a).$

If $x\in B,$ say $x=b\in B,$ let $f(b):=f_{B}(b).$

If $x\in R,$ say $x=x_{1}x_{2}...x_{n}$ where $x_{i}\in A\cup B,n\geq 1,$
let $f(x)=f(x_{1})f(x_{2})...f(x_{n}).$ This is well-defined by the previous
two steps.

If $x\in R^{-},$ say $x=-r$ for some $r\in R,$ let $f(-r):=-f(r);$ the
additive inverse of $f(r)$ in $N$. This does not contradict, for example, $%
-a\in A$ since $f_{A}(-a)=-f_{A}(a)$ as $f_{A}$ is a homomorphism.

If $x\in B_{0}=Fg(R,+),$ say $x=x_{1}+x_{2}+...+x_{n}$ where $x_{i}\in R\cup
R^{-},n\geq 1,$ let $f(x):=f(x_{1})+f(x_{2})+...+f(x_{n})$ which is a
well-defined sum in $N$ from the above.

\noindent (b) Suppose $f(y)\in N$ is defined for all $y\in B_{n-1}$ for some 
$n\geq 1.$ Then we define $f(x)$ for $x\in B_{n}-B_{n-1}=C_{n}\cup D_{n-1}$
considering the two cases:

$(i)$ If $x\in D_{n-1}=R\cdot C_{n-1},$ say $x=r\cdot c=rc$ for some $r\in
R,c\in C_{n-1}\subseteq B_{n-1}$ by the definition of the product in $F(R).$
The induction assumption and part (a) ensure that $f(x)=f(rc):=f(r)f(c)$ is
a well-defined element of $N.$

$(ii)$ If $x\in C_{n}\subseteq B_{n}=B_{n-1}\ast Fg(D_{n-1}),$ then $%
x=x_{1}+x_{2}+...+x_{n}$ where $x_{i}\in B_{n-1}\cup D_{n-1}\cup
D_{n-1}^{-},n\geq 1.$ In this case, let $%
f(x):=f(x_{1})+f(x_{2})+...+f(x_{n}) $ which is well-defined by the
induction assumption and case $(i),$ the latter also justifies the
definition $f(-d):=-f(d)$ for all $d\in D_{n-1}.$

\bigskip

\noindent Thus $f:F(R)\rightarrow N$ is a well-defined function and clearly $%
f\circ \iota _{A}=f_{A}$ and $f\circ \iota _{B}=f_{B}.$ It is
straightforward to check that $f$ is a nearring homomorphism. If $%
g:F(R)\rightarrow N$ is any other nearring homomorphism for which $g\circ
\iota _{A}=f_{A}$ and $g\circ \iota _{B}=f_{B},$ then $g=f$ can easily be
shown since both $f$ and $g$ are homomorphisms. Thus $F(R)$ is the free
product of the nearrings $A$ and $B.$

\bigskip

\noindent We conclude with another example mentioned by Meldrum [5],
extending to nearrings the canonical method used in ring theory for
adjoining an identity. In this example, we construct the free product of a
nearring $A$ with the ring of integers $%
\mathbb{Z}
,$ regarded as a nearring, in which the identity $1\in 
\mathbb{Z}
$ is the identity in the free product by amalgamating both $a1$ and $1a$
with $a$ for all $a\in A.$ Actually, for the ring case, where the canonical
unital extension of a ring $A$ is often called the Dorroh extension of $A$
and is denoted by $D(A),$ we have that $D(A)$ is a semidirect sum of $A$ and 
$%
\mathbb{Z}
$ containing $A$ as an ideal with $D(A)/A\cong 
\mathbb{Z}
$. Such a unital extension of a nearring $A$ containing $A$ as an ideal is
also possible provided the nearring $A$ fulfills a number of requirements.
For this, see Betsch [1] or Veldsman [8]. For the free product in the
example below, there are no restrictions on the nearring $A$ and it is not
necessarily an ideal in the free product.

\bigskip

\noindent \textbf{Example. }For the nearring $A=(A,+,\cdot ,0_{A})$ and the
ring of integers $%
\mathbb{Z}
=(%
\mathbb{Z}
,+,\cdot ,0_{%
\mathbb{Z}
},1_{%
\mathbb{Z}
}),$ let $R$ be the free product of the two semigroups with left zero $%
(A,\cdot ,0_{A})$ and $(%
\mathbb{Z}
,\cdot ,0_{%
\mathbb{Z}
})$ amalgamating $1_{%
\mathbb{Z}
}a$ and $a,$ and $a1_{%
\mathbb{Z}
}$ and $a$ for all $a\in A.$ An explicit construction of $R$ is given here.

\noindent A word $w=x_{1}x_{2}...x_{n}$ over $A\cup 
\mathbb{Z}
,$ $x_{i}\in A\cup 
\mathbb{Z}
,$ is reduced if:

(i) If $x_{i},x_{i+1}\in A$ $($respt. $%
\mathbb{Z}
)$ and $x_{i}x_{i+1}=x_{i}^{\prime }\in A$ $($respt. $%
\mathbb{Z}
),$ replace $x_{i}x_{i+1}$ in $w$ with $x_{i}^{\prime }.$

(ii) If $x_{i}\in \{0_{A},0_{%
\mathbb{Z}
}\},$ replace $x_{1}x_{2}...x_{n}$ with $x_{1}x_{2}...x_{i-1}0_{A}.$

(iii) If $x_{i}=1_{%
\mathbb{Z}
}$ for some $i,$ replace $x_{i-1}x_{i}$ with $x_{i-1}$ in $w$ and replace $%
x_{i}x_{i+1}$ with $x_{i+1}$ in $w.$

\noindent Let $R$ be the set of all reduced words over $A\cup 
\mathbb{Z}
$ with operation concatenation followed by reduction. Then $R$ is a
semigroup with left zero $0:=0_{A}(=0_{%
\mathbb{Z}
})$ and identity $1:=1_{%
\mathbb{Z}
}.$

\noindent Using this $R$ and replacing $D_{n}=R\cdot C_{n}$ with $%
D_{n}:=(R-\{1\})\cdot C_{n}$ for all $n\geq 0,$ proceed as in 2.2 and the
subsequent steps to construct $F(R),$ and using $1$ as the multiplicative
identity throughout the process. Then $F(R)$ is the free product of the
nearrings $A$ and $%
\mathbb{Z}
$ having an identity and containing $A$ as a subnearring.

\bigskip

\noindent \textbf{References}

\noindent \lbrack 1] G. Betsch. Embedding of a Near-Ring into a Near-Ring
with Identity, \textit{North-Holland Mathematics Studies} Volume \textbf{137}
(1987), 37-40.

\noindent \lbrack 2] J.R. Clay. \textit{Nearrings: Geneses and Applications}%
. Oxford University Press, Oxford (1992).

\noindent \lbrack 3] C.C. Ferrero and G. Ferrero. \textit{Nearrings: Some
Developments Linked to Semi-groups and Groups}. Kluwer Academic Publication,
Netherlands (2002).

\noindent \lbrack 4] R.R. Laxton and R. Lockheart. The near-rings hosted by
a class of groups, \textit{Proc. Edinburgh Math. Soc}. \textbf{23} (1980),
69-86.

\noindent \lbrack 5] J.D.P. Meldrum. Free products of near-rings and their
modules, \textit{Algebra Universalis} \textbf{23} (1986), 123-131.

\noindent \lbrack 6] J.D.P. Meldrum. \textit{Near-rings and their links with
groups}. Research Notes in Mathematics 134, Pitman Publishing Limited,
London (1985).

\noindent \lbrack 7] G. Pilz. \textit{Near-rings: The Theory and Its
Applications,} vol. 23. North Holland Mathematics Studies, Amsterdam (1983)

\noindent \lbrack 8] S. Veldsman. On unital extensions of near-rings and
their radicals, \textit{Math. Pannonica} \textbf{3 }(1992), 77-81.

\end{document}